\documentclass[11pt]{article}

\usepackage{makeidx}
\usepackage{pb-diagram}
\usepackage{lamsarrow}
\usepackage{pb-lams}
\usepackage{amsmath}
\usepackage{amsthm}
\usepackage{epsfig}
\usepackage{amssymb}

\theoremstyle{plain}
\newtheorem{theorem}{Theorem}[section]
\newtheorem{lemma}[theorem]{Lemma}
\newtheorem{proposition}[theorem]{Proposition}

\theoremstyle{definition}

\newcommand\bH{{\mathbb H}}

\newcommand{\bQ}{{{\mathbb Q}}}

\newcommand\bZ{{\mathbb Z}}

\renewcommand{\Box}{\blacksquare}

\newcommand{\gF}{{\mathfrak F}}

\newcommand{\Z}{{\mathfrak Z}}

\newcommand{\s}{{\mathcal S}}

\newcommand{\ra}{\rightarrow}
\newcommand{\hra}{\hookrightarrow}
\newcommand{\Lra}{{\longrightarrow}}

\def\inpr{\mathbin{\hbox to 6pt{\vrule height0.4pt width5pt depth0pt \kern-.4pt \vrule height6pt width0.4pt depth0pt\hss}}}

\newcommand{\si}{{\sigma}}

\makeindex

\begin{document}

\title{On the  Reidemeister torsion of rational homology spheres}

\author{Liviu I. Nicolaescu \\University of Notre Dame\\Notre Dame, IN 46556\\
http://www.nd.edu/ $\tilde{}$ lnicolae/ }

\date{June 2000}

\maketitle

\section*{Introduction}

In   the paper \cite{Tu5} V. Turaev  has   proved   a certain identity involving the Reidemeister torsion  of a rational homology   sphere. In this very short note  we will suitably interpret  this identity  as a  second order finite difference equation satisfied by the torsion which will allow  us to prove   a general structure result for the $\mod {\bZ}$ reduction of the torsion.   More precisely  we prove that the mod ${\bf Z}$ reduction of the torsion is completely determined by three data.

\medskip

\noindent $\bullet$ a certain canonical $spin^c$ structure,

\noindent $\bullet$ the linking form of the rational homology sphere and

\noindent $\bullet$ a constant $c\in {\bQ}/{\bZ}$.

\medskip

As a consequence, the constant $c$ is  a  ${\bQ}/{\bZ}$-valued invariant of the rational homology sphere. Experimentations with lens spaces suggest this invariant is as powerful as the torsion itself.

\tableofcontents

\section{The Reidemeister torsion}
\setcounter{equation}{0}

We review briefly a few basic facts about the Reidemeister torsion a  rational homology $3$-sphere. For more details and examples we refer to \cite{Nicol, Tu5}.

 Suppose $M$ is a rational homology sphere. We set $H:= H_1(M, {\bZ})$  and use the multiplicative notation to denote the  group operation on $H$. Denote $Spin^c(M)$ the  $H$-torsor of isomorphism classes of  $spin^c$ structure on $M$.          We denote by  ${\gF}$  the space of functions
\[
\phi: H\ra {\bQ}.
\]
The group $H$ acts on  ${\gF}_M$  by
\[
H\times {\gF}\ni (g,\phi)\mapsto g\cdot \phi
\]
where
\[
(g\cdot\phi)(h)= \phi(hg).
\]
We denote by $\int_H$ the augmentation map
\[
{\gF}_M\ra {\bQ},\;\;\int_H \phi:=\sum_{h\in H}\phi(h).
\]
According to \cite{Tu5} Reidemeister torsion  is a $H$-equivariant map
\[
\tau: Spin^c(M)\ra {\gF}_M,\;\;Spin^c(M)\si\mapsto \tau_\si\in {\gF}_M
\]
such that
\[
\int_H\tau_\si=0
\]
Denote by ${\bf lk}_M$ the linking form of $M$,
\[
{\bf lk}_M: H\times H\ra {\bQ}/{\bZ}.
\]
V. Turaev has proved in \cite{Tu5}  that  $\tau_\si$ satisfies the identity
\begin{equation}
\tau_\si(g_1g_2h)-\tau_\si(g_1h)-\tau_\si(g_2h)+\tau_\si(h)=-{\bf lk}_M(g_1,g_2)\mod {\bZ}
\label{eq: tur}
\end{equation}
$\forall g_1,g_2,h\in H,\;\si\in Spin^c(M)$.

\section{A second order ``differential equation''}
\setcounter{equation}{0}

The identity (\ref{eq: tur}) admits a more suggestive interpretation. To describe it we need a few more notation.

Denote by $\s$ the   space of functions $H\ra {\bQ}/{\bZ}$.   Each $g\in H$ defines a  first order differential operator
\[
\Delta_g:\s\ra\s,\;\;(\Delta_g u)(h):=u(gh)-u(h),\;\;\forall u\in \s,\;\;h\in H.
\]
If $\Xi=\Xi_\si$ denotes the mod ${\bZ}$ reduction of  $\tau_\si$ then  we can rewrite  (\ref{eq: tur}) as \begin{equation}
(\Delta_{g_1}\Delta_{g_2}\Xi)(h)=-{\bf lk}_M(g_1,g_2)
\label{eq: second}
\end{equation}

We will prove  uniqueness and existence results for this equation. We begin with the (almost) uniqueness part.

\begin{lemma}
The  second order linear  differential equation  (\ref{eq: second})  determines $\Xi$ up to an ``affine''  function.
\end{lemma}

\noindent {\bf Proof}\hspace{.3cm} Suppose $\Xi_1$, $\Xi_2$ are two solutions of  the above equation.  Set $\Psi:=\Xi_1-\Xi_2$.   $\Psi$ satisfies the  equation
\[
\Delta_{g_1}\Delta_{g_2}\Psi=0.
\]
Now observe  that any function $F\in \s$ satisfying the second order equation
\[
\Delta_u\Delta_vF=0,\;\;\forall u,v\in H
\]
is affine, i.e. it has the form
\[
F=c + \lambda
\]
where $c\in {\bQ}/{\bH}$ is a constant and $\lambda:H \ra {\bQ}/{\bZ}$ is a character. Indeed, the condition
\[
\Delta_u(\Delta_vF)=0,\;\;\forall u
\]
implies $\Delta_vF$ is a constant depending on $v$, $c(v)$. Thus
\[
F(vh)-F(h)=c(v),\;\;\forall h.
\]
The function $G= F-F(1)$ satisfies the same differential equation
\[
G(vh)-G(h)=c(v)
\]
and the additional condition $G(1)=0$.   If we set $h=1$ in the above equation we deduce
\[
G(v)=c(v).
\]
Hence
\[
G(vh)= G(h)+ G(v),\;\;\forall v, h
\]
 so that $G$ is a character and $F= F(1)+ G$.  Thus, the  differential equation (\ref{eq: second}) determines $\Xi$ up to a constant and a character. $\Box$

 \vspace{1cm}

\begin{lemma}  Suppose $b:H\times H\ra {\bQ}/{\bZ}$ is a nonsigular, symmetric bilinear form on $H$. Then  there exists a quadratic form $q: H\ra {\bQ}/{\bZ}$ such that
\[
\Delta q=b
\]
where
\[
(\Delta q)(uv):=q(uv)-q(u)-q(v).
\]
\end{lemma}

\noindent{\bf Proof}\footnote{We are indebted to Andrew Ranicki for suggesting this approach.}\hspace{.3cm}   Let us briefly recall the terminology in this lemma.   $b$ is nonsingular if the induced map $G\ra G^\sharp$ is an isomorphism. A quadratic map form is a function $q: H\ra {\bQ}/{\bZ}$ such that
\[
q(0)=0,\;\;q(u^k)=k^2q(u),\;\;\forall u\in H, \;k\in {\bZ}
\]
and $\Delta q$ is a bilinear  form.

Suppose $b$ is a nonsingular, symmetric , bilinear  form $H\times H\ra {\bQ}/{\Z}$.  Then, according to \cite[\S 7]{Wall}, $b$ admits a resolution. This  is a  nondegenerate,  symmetric,  bilinear form
\[
B: \Lambda\times \Lambda \ra {\bZ}
\]
on a free abelian group  $\Lambda$ such that,  the induced monomorphism $J_B:\Lambda\ra \Lambda^*$   is a resolution of $H$
\[
0\hra \Lambda\stackrel{J_B}{\Lra}\Lambda^*\stackrel{\pi}{\twoheadrightarrow }H\ra 0
\]
and  $b$ coincides with the induced bilinear form on $\Lambda^*/(J_B\Lambda)$ ($n:=\# H$)
\[
b(\pi(u), \pi(v))= \frac{1}{n^2}B(J_B^{-1}(nu), J_B^{-1}(nv))\mod {\bZ},\;\;\forall u, v\in \Lambda^*.
\]
 Now set
\[
q(\pi(u))=\frac{1}{2n^2}B(J_B^{-1}(nu),J_B^{-1}(nu))\;\;\mod{\bZ}
\]
It is clear that this quantity is well defined i.e.
\[
\frac{1}{2n^2}B(J_B^{-1}(nu),J_B^{-1}(nu))=\frac{1}{2n^2}B(J_B^{-1}(nv),J_B^{-1}(nv))\;\mod{\bZ}
\]
if $v= u+ J_B\lambda, \;\lambda\in \Lambda$. Clearly
\[
\Delta q= b.\;\;\;\Box
\]

\vspace{.7cm}

We deduce that  there exists a constant $c$, a character $\lambda: H\ra {\bQ}/{\bZ}$ and a quadratic form $q$ such that
\[
\Xi(h)=\Xi_\si (h)= c+\lambda(h) + q(h),\;\;\Delta q= {\bf lk}_M.
\]
In the above discussion the choice of the $spin^c$ structure $\si$ is tantamount to a choice of  an origin of $H$ which allowed us to identify the torsion of $M$ as a function $H\ra {\bQ}$. Once we make such a non-canonical choice, we have to replace  $\Xi$  with the family of translates
\[
\{\Xi_g(\bullet):=\Xi(g\bullet);\;\;g\in H\}
\]
In particular
\[
\Xi_g(h):=\Xi(gh)= c+\lambda(gh) + q(gh)= \underbrace{\Bigl(c+ \lambda(g) + q(g)\Bigr)}_{c(g)} + \underbrace{\Bigl(\lambda(h) + (\Delta q)(g,h)\Bigr)}_{\lambda_g(h)} + q(h)
\]
where $\lambda_g(\bullet)= \lambda(\bullet)+ {\bf lk}_M(g, \bullet)$. Since the linking from is nondegenerate  we can find  an {\em unique} $g$ such that $\lambda_g=0$.

We have proved the following  result.

\begin{proposition} Suppose $M$ is a  rational homology sphere. Then there exists an unique $spin^c$-structure $\si$ on $M$ so that, with respect to this  choice  the mod ${\bZ}$ reduction of $\tau_{M,\si}$
\[
\Xi(h):=\tau_{\si}(h)\;\mod {\bZ}
\]
has the form
\[
\Xi(h)= c+ q(h)
\]
where $c\in {\bQ}/{\bZ}$ is a constant  while  $q(u)$ is the unique quadratic form such that
\[
\ \Delta q =-{\bf lk}_M.
\]
In particular,
\[
\Xi(h)=\Xi(h^{-1})\;\mod{\bZ},
\]
and the constant $c\in {\bQ}/{\bZ}$ is a topological invariant of $M$.
\label{prop: new}
\end{proposition}

\section{Examples}
\setcounter{section}{0}

We want to show on  some simple examples that the invariant $c$ is nontrivial.

\medskip

(a) Suppose $M= L(8,3)$. Then its torsion is  (see \cite{Nico2})
 \[
{T_{8, \,3}} \sim  - {\displaystyle \frac {9}{32}} \,x^{7} -
{\displaystyle \frac {3}{32}} \,x^{6} - {\displaystyle \frac {9}{
32}} \,x^{5} + {\displaystyle \frac {5}{32}} \,x^{4} +
{\displaystyle \frac {7}{32}} \,x^{3} - {\displaystyle \frac {3}{
32}} \,x^{2} + {\displaystyle \frac {7}{32}} \,x +
{\displaystyle \frac {5}{32}}
\]
where $x^8=1$ is a generator of ${\bZ}^8$.  Then
\[
q(x^n) = \frac{-3k^2n^2}{16}
\]
The set of  possible  values $\frac{-3m^2}{16}$  mod ${\bZ}$ is
\[
A:=\{0,\;\; \frac{-3}{16},\;\;\frac{4}{16},\;\;\frac{5}{16}\}
\]
The set  possible values of $\Xi(h)$ is
\[
B:=\{-\frac{9}{32},\;\;-\frac{3}{32},\;\;\frac{5}{32},\;\;\frac{7}{32}\}.
\]
We need to find a constant $c\in {\bQ}/{\bZ}$ such that
\[
B-c=A.
\]
Equivalently, we need to figure out orderings $\{a_1, a_2,a_3,a_4\}$    and $\{b_1, b_2, b_3, b_4\}$ of $A$ and $B$ such that $b_i-a_i$ mod ${\bZ}$ is a constant independent of $i$. A little trial and error shows that
\[
\vec{A}=(0, -\frac{3}{16},\frac{4}{16},\frac{5}{16}),\;\; \vec{B}=(-\frac{3}{32},-\frac{9}{32}, \frac{5}{32},\frac{7}{32})
\]
and the constant is $c=-\frac{3}{32}$. This is the coefficient of $x^2$. We deduce that (modulo ${\bZ}$)
\[
F:=T_{8,3}(x)+\frac{3}{32}\sim - {\displaystyle \frac {3}{16}} \,x^{7} -
0\cdot x^{6} - {\displaystyle \frac {3}{16}} \,x^{5} + {\displaystyle \frac {1}{4}} \,x^{4} +
{\displaystyle \frac {1}{4}} \,x^{3} - 0\cdot x^{2} + {\displaystyle \frac {1}{4}} \,x +
{\displaystyle \frac {1}{4}}
\]
The translation of $F$ by $x^{-2}$ is
\[
x^{-2}(T_{8,3}+\frac{3}{32})=\frac{1}{4}x^7+\frac{1}{4}x^6-\frac{3}{16}x^5-\frac{3}{16}x^3 +\frac{1}{4}x^2 +\frac{1}{4}x.
\]

(b) Suppose $M=L(7,2)$. Then, its torsion is (see \cite{Nico2})
\[
{T_{7, \,2}} \sim  - {\displaystyle \frac {2}{7}} \,x^{6} +
{\displaystyle \frac {1}{7}} \,x^{5} + {\displaystyle \frac {2}{7
}} \,x^{3} + {\displaystyle \frac {1}{7}} \,x - {\displaystyle
\frac {2}{7}}
\]
where $x^7=1$ is a generator of ${\bZ}_7$. We see that  in this form $T_{7,2}$ is symmetric, i.e. the coefficient of $x^k$ is equal to the coefficient of $x^{6-k}$. The constant $c$ in this case must be the   coefficient of the middle  monomial $x^3$,  which is $\frac{2}{7}$.

(c) Suppose $M=L(7,1)$. Then
\[
{T_{7, \,1}} \sim {\displaystyle \frac {2}{7}} \,x^{6} +
{\displaystyle \frac {1}{7}} \,x^{5} - {\displaystyle \frac {1}{7
}} \,x^{4} - {\displaystyle \frac {4}{7}} \,x^{3} -
{\displaystyle \frac {1}{7}} \,x^{2} + {\displaystyle \frac {1}{7
}} \,x + {\displaystyle \frac {2}{7}}
\]
This is again a symmetric polynomial and   the coefficient of the middle monomial is $-4/7$.  We see that this invariant distinguishes the lens spaces $L(7,1)$, $L(7,2)$.

(d) For $M= L(9,2)$ we have
\[
 {T_{9, \,2}} \sim  - {\displaystyle \frac {10}{27}} \,x^{8} +
{\displaystyle \frac {2}{27}} \,x^{7} - {\displaystyle \frac {1}{
27}} \,x^{6} + {\displaystyle \frac {8}{27}} \,x^{5} +
{\displaystyle \frac {2}{27}} \,x^{4} + {\displaystyle \frac {8}{
27}} \,x^{3} - {\displaystyle \frac {1}{27}} \,x^{2} +
{\displaystyle \frac {2}{27}} \,x - {\displaystyle \frac {10}{27}}
\]
Again, this is a symmetric function, i.e the coefficient of $x^k$  is equal to the coefficient of $x^{8-k}$, $x^9=1$. The constant  is the coefficient of $x^5$, which is $2/27$. We deduce that, mod ${\bZ}$, we have
\[
T_{9,2}=-\frac{2}{3}x^8-\frac{2}{9}x^7-\frac{1}{3}x^6 -\frac{2}{9}x^7
\]

(e) Finally  when $M= L(9,7)$ we have
\[
{T_{9, \,7}} \sim  - {\displaystyle \frac {8}{27}} \,x^{8} -
{\displaystyle \frac {2}{27}} \,x^{7} + {\displaystyle \frac {10
}{27}} \,x^{6} + {\displaystyle \frac {1}{27}} \,x^{5} -
{\displaystyle \frac {2}{27}} \,x^{4} + {\displaystyle \frac {1}{
27}} \,x^{3} + {\displaystyle \frac {10}{27}} \,x^{2} -
{\displaystyle \frac {2}{27}} \,x - {\displaystyle \frac {8}{27}
}
\]
the  polynomial is again symmetric so that the constant $c$ is the coefficient of $x^4$ which is $-2/7$.

\bigskip

It would be very interesting to know whether the invariant $c$ satisfies  any surgery properties.  This is not a trivial issue because     we cannot  relate  the potential surgery properties of $c$ to the surgery properties of the torsion.  In the case of torsion the surgery formula involve finite difference operators which  kill the constants so $c$ will not appear in any of them.

\end{document}